\documentclass[12pt]{amsart}
\usepackage{amscd,amsmath,amssymb,amsfonts}
\theoremstyle{plain}
\newtheorem{thm}{Theorem}

\theoremstyle{definition}

\newtheorem{rmk}[thm]{Remark}

\numberwithin{thm}{section}
\numberwithin{equation}{section}

\newcommand{\ml}[2]{\begin{multline}\label{#1}#2 \end{multline}}
\newcommand{\ga}[2]{\begin{gather}\label{#1}#2 \end{gather}}

\newcommand{\sL}{{\mathcal L}}

\newcommand{\sO}{{\mathcal O}}
\newcommand{\sP}{{\mathcal P}}

\newcommand{\C}{{\mathbb C}}

\newcommand{\F}{{\mathbb F}}

\begin{document}

\title[Divisibility]{Some elementary theorems about divisibility of 0-cycles 
on abelian varieties defined over finite fields}
\author{H\'el\`ene Esnault}
\address{Mathematik,
Universit\"at Essen, FB6, Mathematik, 45117 Essen, Germany}
\email{esnault@uni-essen.de}
\date{Oct. 23, 2003}
\begin{abstract}
If $X$ is an abelian variety over a field, and $L\in {\rm Pic}(X)$, we
  know that the degree of the 0-cycle $L^g$ is divisible by $g!$. 
As a 0-cycle, it is not, even over a field of cohomological dimension
  1, but we show that over a finite field, there is some hope.

\end{abstract}
\maketitle
\begin{quote}

\end{quote}

\section{Introduction}
Let $X$ be an abelian variety defined over a field $k$, and let $L\in
{\rm Pic}(X)$. Then the Riemann-Roch theorem  implies that the
degree of 
0-cycle 
$L^g$ is divisible by $g!$. An important example is when $X$ is the
product of an abelian variety with its dual $A\times A^\vee$, $A$ of
dimension $n$, thus $g=2n$, and $L$ is the Poincar\'e
bundle, normalized so as to be trivial along the 2 zero sections. In
this case, the highest product $L^g$ is just $g!$ times
the origin, as a 0-cycle
in
$CH_0(A\times A^\vee)$.  In this note we investigate 
the question of $g!$-divisibility of the 0-cycle $L^g$ in
$CH_0(X)$. We show that already for $g=2$, there are counter-examples
(see Remark \ref{rmk4.1}). 
This answers negatively a question by B. Kahn: there are no divided
powers in the Chow groups of abelian varieties, not even in their
\'etale
motivic cohomology. 
Indeed, his
question was the motivation to study divisibility of the 0-cycle
$L^g$. Our method consists in relating this divisibility for $g=2$ to
the existence of theta characteristics on smooth projective curves
over the given field. Over finite fields, Serre's theorem asserts that
there are always theta characteristics. We derive from this that if
$g=2n$ and the field is finite, 
one always has 2-divisibility (see Theorem \ref{2div}, Remark \ref{rmk4.3}). 
Thus, the 2-divisibility is an arithmetic statement (see Remark
\ref{rmk4.2}). 
In order to find a non-trivial class of invertible sheaves
$L$ for which one has 
$g!$-divisibility, one needs more arithmetic. In Theorem \ref{thm3.1} we
show that a principal polarization $L$ of geometric origin 
has the strong property that $L^g$ is $g!$-divisible as a 0-cycle. 
Aside of Serre's theorem mentioned above, 
the proof relies on (an adaption of) Mattuck's results on geometric
polarization (\cite{M}), on abelian class field theory by
K. Kato-S. Saito
(\cite{KS}) and on S. Bloch's theorem on 0-cycles on abelian
varieties (\cite{B}). Perhaps it gives some hope that $g!$-divisibility is true
in general over a finitie field.

{\it Acknowledgement}: The note is an answer to a question by
B. Kahn. S. Bloch explained
to us long ago the use of abelian class field theory to compute
0-cycles over finite field. This note relies on this idea. 
J.-P. Serre explained to us his theta theorem. We thank them all, and
also E. Howe, A. Polishchuk, and E. Viehweg for answering our
questions. In particular, Remark \ref{rmk4.2} comes from a discussion
with E. Viehweg.

\section{$2$-divisibility}
\begin{thm}\label{2div} Let $X$ be an abelian variety of dimension
  $g=2n$
defined over a finite field. 
Let $L\in {\rm Pic}(X)$ be a line bundle. Then the 0-cycle
$L^g$ is divisible by $2$ in the Chow group of 0-cycles of $X$. \end{thm}
\begin{proof}
If $A$ is a sufficiently
  ample line bundle, then $L\otimes A=B$ is  very ample as well. One has
$L^g\equiv A^g +  B^g \ {\rm mod}  \ 2CH_0(X)$, 
thus one may assume that $L$ is as ample as necessary.
For $L$ sufficiently ample, then there is a finite field extension
$k'\supset k$ of odd degree, such that the intersection
 $C\subset X\times_k k'$  of $(g-1)$ linear sections of $L\times_k k'$ 
in general position is smooth.  One has
$\omega_C=(L\times_k k')^{\otimes (g-1)}|_C \equiv (L\otimes_k k')|_C \ {\rm mod} \
 2 {\rm Pic}(C)$, 
thus via the Gysin homomorphism
${\rm Pic}(C)\xrightarrow{\iota} CH_0(X\times_k k')$
one has
$(L\times_k k')^g\equiv \iota_*(\omega_C) \ {\rm mod} \ 2CH_0(X \times_k k').$
On the other hand, Serre's theorem \cite{S}, p.61, Remark 2,  asserts
  that a smooth curve over a finite field admits a theta
  characteristic, that is that $\omega_C$ is 2-divisible in ${\rm
  Pic}(C)$. This shows that $(L\times_k k')^g \in 2CH_0(X\times_k k')$,
  thus by projection formula, since $k'\supset k$ is odd, that $L^g
\in 2CH_0(X)$ as well. 
\end{proof}
\section{$g!$-divisibility}
Let $C$ be a smooth projective curve of genus $g$ 
defined over a finite field $k=\F_q$.  
Then $C$ has carries a 0-cycle $p$  of degree 1. Let $J$ be the
Jacobian of $C$. The rest of this section is trivial if $g\le 1$ thus
we assume that $g\ge 2$.  
We consider the cycle map
\ga{}{\psi_p: C\to  J, \ y\mapsto \sO_C(y-{\rm deg}(y)\cdot p).\notag}
The cycle map induces a birational morphism which we still denote by $\psi_p$
\ga{}{\psi_p: {\rm Sym}^g(C)\to J, (x_1,\ldots,x_g)\mapsto 
\otimes_{i=1}^g\psi_p(x_i).\notag}
In particular, writing
$p=\sum_i m_i p_i =q_1-q_2$,  with  
$q_i$  effective, one has 
${\rm deg}(q_1)-{\rm deg}(q_2)= \sum_i m_i{\rm deg}(p_i)=1.$
We denote by $ \pi: C^g\to {\rm Sym}^g(C)$ the quotient map. It defines
the divisor
\ga{}{D_p=\sum_i m_i (\psi_p\circ \pi)_*
(C^{g-1}\times p_i), \ L_p=\sO_J(D_p) \in {\rm Pic}(J).\notag}
We know that $L_p$ is a principal polarization, thus the 0-cycle
$L_p^g$ has degree $g!$. 
The purpose of this section is to show
\begin{thm} \label{thm3.1}
The 0-cycle $L_p^g$ in the Chow group $CH_0(J)$ 
of 0-cycles of $J$ is divisible by $g!$ in $CH_0(J)$, 
that is there is a 0-cycle $\xi\in CH_0(J)$ of degree 1 with 
$L^g_p=g! \cdot\xi\in CH_0(J)$. 
\end{thm}
\begin{proof} Since the 0-cycle $p$ of degree 1 won't change during
  the proof, we simplify the notation and set $\psi=\psi_p, D=D_p, L=L_p$. 
We consider the Poincar\'e bundle 
$P=p_1^*L\otimes p_2^*L \otimes \mu^*L^{-1}\in {\rm Pic}(J\times J)$
where 
$\mu: J\times J\to J, \ \mu(x,y)=x+y.$ 
Via the cycle map
\ga{}{ \iota_*:{\rm Alb}(C)={\rm Pic}(C)^0={\rm Pic}(J)^0\xrightarrow{\psi_*}CH_0(J)^{{\rm deg}=0}\xrightarrow{h} J(k),\notag}
where $h$ is the Albanese mapping of $J$,  
one defines
\ga{}{\iota_*\omega_C(-2(g-1)p)=y\in J(k).\notag}
We now adapt \cite{M}, section 6 to the situation where $p$ is not
necessarily
a $k$-rational point of $C$, but only a 0-cycle of degree 1. 
One defines the involution 
\ga{}{\delta: J\to J \notag \\
 x\mapsto -x+y=\tau_y\circ (-1)^*(x)=(-1)^*\circ \tau_{-y}\notag }
where $\tau_y$ is the tranlastion by $y$ while $(-1)$ is the multiplication by $-1$. 
We set
\ga{}{\ell=\psi^*L \notag \\ \sP=(\psi \times 1)^*P=
p_1^*\ell \otimes p_2^*L \otimes (\psi \times 1)^*\mu^*L^{-1}.\notag }
As in \cite{P}, p. 249, for $d>2g-2$, the Riemann-Roch theorem asserts that 
\ga{}{E_d:=p_{2*}(\sP\otimes p_1^*\sO_C(d p))=
Rp_{2*}(\sP\otimes p_1^*\sO_C(d p)) \notag }
is a vector bundle. One has
\ga{}{\otimes_i \sP|_{p_i\times J}^{\otimes m_i}=
\sP|_{C\times \{0\}}=\sO_J. 
\notag}
This implies that for any natural number $f>0$ one has
\ga{3.1}{ \sP_{fq_1\times J}\otimes 
\sP_{fq_2 \times J}^{-1}=\sP|_{C\times \{0\}}=\sO_J. }
For two natural numbers $e> f> 0$, one has the diagram
\ga{}{\begin{CD}
\sP\otimes p_1^*\sO_C((e-f)q_1-eq_2)@> fq_1>> \sP\otimes p_1^*\sO_C(ep)\\
@V fq_2 VV \\
\sP\otimes p_1^*\sO_C((e-f)p)
\end{CD}\notag
}
where the horizontal map is induced by
$\sO_C\to \sO_C(fq_1)$ and the vertical one
by  $\sO_C\to \sO_C(fq_2)$. 
Thus one has the relation in $K_0(J)$
\ga{}{Rp_{2*}( \sP\otimes p_1^*\sO_C(ep))-Rp_{2*}(
\sP\otimes p_1^*\sO_C((e-f)p))= \notag \\
p_{2*}(\sP \otimes p_1^*\sO_C((e-f)p)|_{fq_2 \times J})^{-1}
\otimes p_{2*}(\sP\otimes p_1^*\sO_C(ep)|_{fq_1 \times J})\notag\\
\xrightarrow{\eqref{3.1}\cong}0.
\notag}
For $d>2g-2$, we then have in the Grothendieck group $K_0(J)$ 
\ga{3.2}{p_{2*}(\sP\otimes p_1^*\sO_C(dp))+R^1p_{2*}
(\sP\otimes \sO_C((2g-2-d)p)=0.}
On the other hand, one has
\ga{}{\delta^*E_d=p_{2*}(\sP^{-1}\otimes
  p_1^*\omega_C((d-2g+2)p))\notag } and 
Serre duality implies
\ga{3.3}{\delta^*(E_d^\vee)= R^1p_{2*} (
\sP \otimes p_1^*(\sO((2g-2-d)p))) }
(see \cite{P}, p. 249). 
The involution $\delta$ acts on $J$,
 thus on the Chow groups of $J$.
Thus \eqref{3.2} and \eqref{3.3}
imply
\ga{}{\delta^*c_1(E_d)=c_1(E_d)=:M, \ M^2=\delta^*M^2=
c_2(E_d)+\delta^*c_2(E_d).\notag}
Thus we have
\ga{3.4}{M^g=M^2\cdot M^{g-2}=\\
M^{g-2}c_2(E_d)+ \delta^*M^{g-2}\cdot \delta^*
c_2(E_d)=Z+\delta^*Z\notag}
for 
\ga{}{
Z:=c_2(E_d)\cdot M^{g-2} \in CH_0(J), \ {\rm deg}Z=n, 2n={\rm deg}M^g.\notag}
Let us set
\ga{}{Z'=Z-n \{0\} \in CH_0(J)^{{\rm deg}=0}.\notag}
By abelian class field theory  for 0-cycles on varieties defined over finite fields, \cite{KS}, Corollary p. 274, the Albanese mapping $h$ is an isomorphism. This allows to identify explicitely $\delta^*$ on $CH_0(J)$. 
Indeed, if $W\in CH_0(J)$ has degree $n$, and if $z\in J(k)$, then 
\ml{3.5}{h(\tau_z^*(W)-n\{0\})= \tau_z^*h(W-n\{0\}) + h(\tau_z^*
n \{0\}-n \{0\})=\\
\tau_z^*h(W-n\{0\}) + nz
.}
On the other hand, S. Bloch's theorem \cite{B}, Theorem 3.1 asserts that the 
second Pontryagin product dies in $CH_0(J\times_k \bar{k})$, where $\bar{k}$ is the algebraic closure of $k$, thus by class field theory again \cite{KS}, p. 274 Proposition 9, it dies in $CH_0(X)$. Thus
\ga{3.6}{\tau_z^*a=a \ {\rm for \ all } \ a\in CH_0(X)^{{\rm deg}=0}.}
Thus \eqref{3.5} and \eqref{3.6} imply
\ga{3.7}{h(\tau_z^*(W)-n\{0\})=h(W-n\{0\}) +nz.}
On the other hand
\ga{3.8}{h((-1)^*W-n\{0\})=-h(W-n\{0\}).}
Thus, \eqref{3.4},  \eqref{3.7} and \eqref{3.8} imply 
\ga{3.9}{h(M^g-2n\{0\})= -ny. }
We now apply again Serre's theorem \cite{S}, p. 61, Remark 2,
which asserts that over a finite
field, a smooth projective curve has a theta divisor. Thus $
\omega_C(-(2g-2)p)\in 2{\rm Pic}^0(C)$ and a fortiori via 
the Gysin homomorphism
$\psi_*$ one has  
\ga{3.10}{y=2\xi_0 \in CH_0(J)^{{\rm deg}= 0},  \ {\rm for \ some} \ \xi_0 \in CH_0(J)^{{\rm deg}=0}.}
Thus \eqref{3.9}, \eqref{3.10}, and \cite{KS}, loc. cit. imply
\ga{3.11}{M^g=2n (\{0\}-\xi_0).}
It remains to compare $M^\vee$ and $L$. As a divisor, one has
$M^\vee \otimes_k \bar{k}=\{\sL\in {\rm Pic}^0(J)(\bar{k}), 
\Gamma(C\times_k \bar{k}, \sL((g-1)p))\neq 0\}$. Thus 
$M^\vee \otimes_k \bar{k}=\sO_{J\times_k \bar{k}}(D)$, as we know,
both underlying divisors are physically the same and they are both
reduced. Thus $M^\vee = \sO_J(D)\otimes \sL$ for some $\sL\in {\rm
  Pic}^0(J)$
which is torsion. On the other hand, the map $J\to {\rm Pic}^0(J),
a\mapsto \tau_a^*\sO_J(D)\otimes \sO_J(D)^{-1}$ is an isomorphism over
$\bar{k}$, thus is an isomorphism over $k$. This implies that $M^\vee=\tau_a^*
\sO_J(D)$ for some $a\in J$. Thus \eqref{3.11} implies
\ga{3.12}{L^g=2n \tau_a^*(\xi_0-\{0\}), \  2n=g!. }
This finishes the proof. 

\end{proof}

\section{Remarks}
\begin{rmk} \label{rmk4.1} B. Kahn asked whether \'etale
motivic cohomology of
  abelian varieties has divided powers. Theorem \ref{2div}, read
  backwards, yields a negative answer. 
Let $k$ be a field, and let $C$ be a genus 2 curve defined over this
  field with the two properties that it carries a 0-cycle  $p$
of degree 1 
  and it does not have a theta characteristic. Let $J$ be the Jacobian
of $C$ and $\psi_p=:\psi: C\to J$ be the cycle map assigned to the
  choice of $p$.  
The composite map 
$\iota_*:{\rm Alb}(C)={\rm Pic}(C)^0={\rm
  Pic}(J)^0\xrightarrow{\psi_*}CH_0(J)^{{\rm deg}=0}\xrightarrow{h}
  J(k)$ being an isomorphism, one has $\iota_*\omega_C(2(g-1)p)\in
  J(k)$ is not 2-divisible. Thus a fortiori, the class of the Gysin
  image of $\omega_C$
won't be 2-divisible in any cohomology which has the property that it
  maps to \'etale cohomology and the kernel maps to the Albanese. For
  example, \'etale motivic cohomology.

It remains to give a concrete curve. 
One could take for $k$ the function field of the fine moduli space of
pointed genus 2 curves with some level over a given algebraically
clsoed field $F$. This has transcendence degree 3 over $F$ and is of
course very large. Here is an example due to J.-P. Serre over the
field
$k=\C(t)$ 
of cohomological dimension 1: $C$ is defined by its hyperelliptic
equation
$y^2=x^6-x-t$. It has 2 rational points at $\infty$. The Galois
group of $x^6-x-t$ is the symmetric group in 6 letters, which acts with
2 orbits on the space of theta characteristics over $\bar{k}$, one
with 6 elements and one with 10. 
\end{rmk}
\begin{rmk} \label{rmk4.2} This remark  arose in discussions
  with E. Viehweg in view of Theorem \ref{2div}. If $X$ is a product of curves
$X=C_1 \times \ldots \times C_g$ over a field $k$, then the Pic
  functor is quadratic after Mumford, which means that a line bundle $L$
  on $X$ is a sum of line bundles $L_{ij}$ coming via pull-back
from only 2 factors
  $(ij), i\neq j$.  Thus the expansion of $L^g$ will have two kinds of 
summands. Those of the type $L_{i_1j_1}\cdots L_{i_gj_g}$ with all
pairs $(i_c, j_c)$ being different. The coefficient of such a summand 
is $g!$ thus this term
is $g!$ divisible. 
Then those of the type $L_{i_1j_1}^2\cdots L_{i_aj_a}^2\cdot
L_{i_{a+1}j_{a+1}}
\cdots L_{i_{a+b}j_{a+b}}$ with all pairs $(i_c,j_c)$ being different and 
$2a+b=g$. The coefficient of such a summand is $\frac{g!}{2^a}$. 
Thus $g!$ divisibility of any $L^g$ on a product of $g$ curves 
splits up into two kinds of divisibility. Over any field $k$, 
geometry always forces
$\frac{g!}{2^a}$ divisibility for $a=[\frac{g}{2}]$, where $[c]$ means
the integral part of a real number $c$. On the other
hand, let $k$ be a field which has the property that any curve has a
theta characteristic, e.g. a finite field (\cite{S}, p.61, Remark
2). Then the argument of Theorem \ref{2div} implies that if $L\in {\rm
  Pic}(C_1\times C_2)$, then $L^2$ is 2-divisible. Indeed, one reduces
as in the proof of Theorem \ref{2div} to the case where $L=\sO(\Gamma)$ 
for a smooth curve $\Gamma \subset C_1\times C_2$ and by the
adjunction
formula 
one has $L^2=i_{\Gamma *}\omega_\Gamma - p_1^*\omega_{C_1}\cdot \Gamma
 -p_2^*\omega_{C_2}\cdot \Gamma $ where $p_i: C_1\times C_2\to C_i$
 and $i_\Gamma: \Gamma \to C_1\times C_2$ is the closed embedding. 
Thus over such a field, $L^g$ is always $g!$-divisible in $CH_0(C_1
\times \ldots \times C_g)$. 

\end{rmk}

\begin{rmk} \label{rmk4.3}
The conclusion of Theorem \ref{2div} is of course true over 
any field $k$ over which any smooth projective curve has a theta
characteristic. 

\end{rmk}

\newpage
\bibliographystyle{plain}
\renewcommand\refname{References}

\end{document}